\documentclass[oneside,a4paper,12pt,reqno]{amsart}
\usepackage{amssymb,amsthm,amsmath,amsfonts,array,amscd}
\usepackage[numbers]{natbib}
\usepackage[colorlinks,citecolor=blue,urlcolor=blue]{hyperref}
\usepackage{hypernat}
\usepackage{geometry}
\usepackage{graphics}
\usepackage[dvips]{graphicx}
\begin{document}

\newtheorem{theorem}{Theorem}
\newtheorem{lemma}{Lemma}
\theoremstyle{remark}
\newtheorem{rem}{Remark}

\def\MR#1{\href{http://www.ams.org/mathscinet-getitem?mr=#1}{MR#1}}

\title{ $U-$max-Statistics and Limit Theorems for \\ Perimeters  and Areas of Random Polygons}
\author{E.~V.~KOROLEVA, \,  Ya.~Yu.~NIKITIN }

\medskip

\noindent \email{ evgeniya.koroleva@hotmail.com, yanikit47@gmail.com }
\date{}
\begin{abstract}
Recently W. Lao  and M. Mayer \cite{lao},  \cite{diss_Lao}, \cite{diss_Mayer}
considered $U$-max - statistics,
where instead of sum  appears the {\it maximum} over the same set of
indices. Such statistics often appear  in stochastic geometry. The
examples are given by the largest distance between
random points in a ball, the maximal diameter
of a random polygon, the largest scalar product within a sample of points, etc.
Their limit distribution is  related to the distribution of extreme values.

Among the interesting results obtained in \cite{lao}, \cite{diss_Lao},
\cite{diss_Mayer} are limit theorems for the maximal
perimeter and the maximal area of random {\it triangles} inscribed in a
circumference. In the present paper we generalize these theorems to
convex $m$-{\it polygons}, $m \geq 3,$ with random vertices on the circumference. Next, a similar
problem is solved for the minimal perimeter and the minimal area
of \textit{circumscribed} $m$-polygons which has not been previously considered in literature.
Finally, we discuss the obtained results when $m \to \infty.$
\end{abstract}



\maketitle


\section{Introduction}

\subsection{One-sample $U$-statistics and parametric functionals}

Let $\xi_1,\ldots,\xi_n$ be inde\-pen\-dent random
elements taking values in a measurable space $(\mathfrak{X}, \mathfrak{A})$ and
having identical distribu\-ti\-on $P.$ Let $\mathcal{P}=\{P\}$ be some class of probability distributions on
$(\mathfrak{X}, \mathfrak{A})$ and let $\theta(P)$ be a functional on $\mathcal{P}$.

The functional $\theta(P)$ is called \textit{regular} \cite{koroluk}, if
 $\theta(P)$ can be represented as
 \begin{equation}
 \label{theta}
\theta(P)=\int\limits_{\mathfrak{X}}^{} \ldots\int\limits_{\mathfrak{X}}^{} h(x_1, \ldots, x_m )P(dx_1)\ldots P(dx_m)
\end{equation}
with some real-valued symmetric Borel function  $h(x_1, \ldots, x_m)$ which is called the kernel, while the
integer number $m\geq 1$ is called the degree of functional $\theta(P)$.

 Halmos and Hoeffding \cite{Hal}, \cite{Hoe}  began to study the class of unbiased
estimates of $\theta(P)$ called $U$-statistics, which are defined as  follows.
Consider a kernel
$h(x_1,\ldots,x_m)$ of parametric functional (\ref{theta}). Then the
{\it $U$-statistic of degree $m$} is defined as
$$U_n={\binom{n}{m}}^{-1}\sum\limits_{J}h(\xi_{i_1}, \ldots, \xi_{i_m} ),$$
where $n\geq m$ and $J=\{(i_1, \ldots,i_m): 1\leq i_1< \ldots< i_m\leq n \}$ is
a set of increasing permutations of indices $i_1, \ldots, i_m.$

It turns out  that numerous statistical estimates
and test statistics belong to the class of $U$-statistics. This entailed
the intensive development of the theory,  see
\cite{koroluk} and \cite{Lee}.

\subsection{$U$-max-statistics}

$U$-max-statistics appear  for the description of extreme counterparts of $U$-statistics,
they are given by the formula:
$$H_n=\max\limits_{J} h(\xi_{i_1},\ldots,  \xi_{i_m}).$$
$U$-min-statistics are
 defined similarly, they can be reduced
to $U$-max-statistics by changing the sign of the kernel.
Here are some examples of $U$-max and $U$-min-statistics.

\smallskip

1. Largest interpoint distance $\max\limits_{1\leq i<j\leq n}^{}\|\xi_i-\xi_j\|,$
\noindent where $\xi_1, \xi_2, \ldots$ are i.i.d. points in  the $d$-dimen\-sio\-nal unit ball $\mathbb{B}^d, d\geq 2$.

2.  Largest scalar product $\max\limits_{1\leq i<j\leq n}^{}\langle \xi_i, \xi_j\rangle,$
where $\xi_1, \xi_2, \ldots$ are i.i.d. points in the ball $\mathbb{B}^d, \  d\geq 2$.

3. Smallest spherical distance:
$\min\limits_{1\leq i<j\leq n}^{}\beta_{i,j},$
where $\beta_{i,j}$ denotes the smaller of two central angles
between $U_i$ and $U_j.$ It is assumed that the  vertices $U_1,\ldots,
U_n$ are i.i.d. points on the unit sphere
$\mathbb{S}^{d-1}, d\geq 2$.

4. Largest perimeter
$\max\limits_{1\leq i<j<l\leq n}^{} peri(U_i, U_j,  U_l)$
and largest area
$\max\limits_{1\leq i<j<l\leq n}^{} area(U_i, U_j,  U_l)$
among all  inscribed triangles, whose vertices are formed by
triplets of  points taken  from a sample $U_1,  \ldots,  U_n$ of indepen\-dent and uniformly distributed points on the unit
circum\-fe\-rence.

Lao and Mayer   \cite{lao}, \cite{diss_Lao}, \cite{diss_Mayer} initiated the study of $U$-max-statistics
and proved for them the basic limit theorem. They used some modification of a statement on Poisson convergence from the monograph of Barbour, Holst and Janson \cite[p.35]{Barb}.

\smallskip

\textbf{Lao-Mayer Theorem \cite{lao}}. \textit{Let $\xi_1,\ldots,
\xi_n$ be i.i.d. random elements in some measurable space
$({\mathfrak X}, \mathfrak {A} )$ and let $h$ be a symmetric Borel function,
$h:{\mathfrak X}^m \rightarrow \mathbb{R}$. Put
$$H_n=\max\limits_{J} h(\xi_{i_1},\ldots,  \xi_{i_m}),$$
and denote for any $z\in\mathbb{R}$}
$$p_{n,z}=P\{h(\xi_1,\ldots,
\xi_m)>z\}, \qquad  \lambda_{n,z}=\binom{n}{m}p_{n,z},$$
$$\tau_{n,z}(r)=P\{h(\xi_1,\ldots,  \xi_m)>z, \, h(\xi_{1+m-r},
\xi_{2+m-r},\ldots,  \xi_{2m-r})>z\}/p_{n,z}.$$

\medskip

\noindent \textit{Then for
all $n\geq m$ and for each $z\in\mathbb{R}$ we have:}
$$|P\{H_n\leq z\}-\exp\{-\lambda_{n,z}\}|\leq$$
\begin{equation}
\label{Lao_Mayer_Th}
\leq(1-\exp\{-\lambda_{n,z}\})\left\{p_{n,z}\left[\binom{n}{m}-\binom{n-m}{m}\right]
+\sum_{r=1}^{m-1}\binom{m}{r}\binom{n-m}{m-r}\tau_{n,z}(r)\right\}.
\end{equation}

Clearly, the result can be reformulated for the minimal value of
the kernel by replacing $h$ with $-h$.

\begin{rem}
\label{remark_1} \textit{(Lao-Mayer \cite{lao}).} If the sample
size $n$ tends to infinity, then the error in (\ref{Lao_Mayer_Th}) is
asymptotically of order
$$O\left(p_{n,z}n^{m-1}+\sum\limits_{r=1}^{m-1}\tau_{n,z}(r)n^{m-r}\right),$$
where for $m>1$ the sum is dominating, see \cite{Barb}.
\end{rem}

Silverman and Brown \cite{Silw} formulated the
conditions which ensure that the general theorem
used in \cite{lao} provides a non-trivial Weibull limit law.

\smallskip

\textbf{Silverman-Brown Theorem \cite{Silw}}. \textit{In the
setting of Lao-Mayer theorem, if for some sequence of transformations $z_n:
T\rightarrow\mathbb{R}, \ T\subset\mathbb{R},$  the conditions:
\begin{equation}
\label{corollary_silv_1}
 \lim\limits_{n\rightarrow\infty}\lambda_{n,z_n(t)}=\lambda_t>0,
\end{equation}
\begin{equation}
\label{corollary_silv_2}
\lim\limits_{n\rightarrow\infty}n^{2m-1}p_{n,z_n(t)}\tau_{n,z_n(t)}(m-1)=0,
\end{equation}
hold  for each $t\in T$, then
\begin{equation}
\label{corollary_silv_3}
\lim\limits_{n\rightarrow\infty}P\{H_n\leq
z_n(t)\}=\exp\{-\lambda_t\}
\end{equation}
for each $t\in T.$}

\begin{rem}
\label{remark_2} \textit{(Lao-Mayer \cite{lao}).}  If $m>2$, then
the condition (\ref{corollary_silv_2}) can be replaced by the
weaker condition:
\begin{equation}
\label{corollary_silv_4}
\lim\limits_{n\rightarrow\infty}n^{2m-r}p_{n,z_n(t)}\tau_{n,z_n(t)}(r)=0
\end{equation}
for all $r\in\{1, \ldots, m-1\}$.
\end{rem}

\begin{rem}
\label{remark_3} \textit{(Lao-Mayer \cite{lao}).} Condition
(\ref{corollary_silv_1}) implies $p_{n,z}=O(n^{-m})$, and
therefore (\ref{corollary_silv_3}) is valid with the rate of convergence
$$
O\left(n^{-1}+\sum\limits_{r=1}^{m-1}n^{2m-r}p_{n,z_n(t)}\tau_{n,z_n(t)}(r)\right).
$$
\end{rem}

For the perimeter and the area of inscribed triangles (see example 4)
Lao and Mayer in \cite{lao}, \cite{diss_Lao} and \cite{diss_Mayer}  obtained
the following results.

\renewcommand{\thetheorem}{\Alph{theorem}}
\renewcommand{\thelemma}{\Alph{lemma}}

\begin{theorem}[Perimeter of inscribed triangle]
\label{theorem_A} Let $U_1, U_2, \ldots$ be independent and
uniformly distributed points on the unit circumference
$\mathbb{S},$ and let $peri(U_i,U_j,U_l)$ be the perimeter of triangle
formed by the triplet of points  $U_i,U_j,U_l.$  Set
$H_n=\max\limits_{1\leq i<j<l\leq n}peri(U_i,U_j,U_l).$  Then for
each $t>0$
$$\lim\limits_{n\rightarrow\infty}P\{n^3(3\sqrt{3}-H_n)\leq
t\}=1-\exp\left\{-\frac{2t}{9\pi}\right\}.$$ The rate of
convergence is $O\left(n^{-\frac{1}{2}}\right)$.
\end{theorem}

As a comment to this result, we note that among all triangles
inscribed in the unit circumference, the
 regular triangle has the
maximal value of perimeter equal to $3\sqrt{3}.$ It is a classical
and well-known result, see \cite{yaglom}. Clearly, the maximal perimeter
$H_n$ of random triangle tends to this value. The theorem gives the
required normalization for this convergence, describes the limit
distribution and establishes the rate of convergence.

Theorem \ref{theorem_A} is proved by means of the following Lemma
\ref{lemma_A}.

\begin{lemma}
\label{lemma_A} Let $U_1, U_2, U_3$ be independent and uniformly
distributed points on the unit circumference $\mathbb{S}.$ Then
$$\lim\limits_{s\rightarrow +0} s^{-1}P\{peri(U_1,U_2,U_3)\geq
3\sqrt{3}-s\}=\frac{4}{3\pi}.$$
\end{lemma}

Now we proceed to random areas.

\begin{theorem}[Area of inscribed triangle]
\label{theorem_B}
Let $U_1, U_2, \ldots$ be independent and
uniformly distributed points on the unit circumference $\mathbb{S}.$ Set
$G_n=\max\limits_{1\leq i<j<l\leq n}area(U_i,U_j,U_l),$ where
$area(U_i,U_j,U_l)$ is the area of random triangle formed by the triplet of points $U_i,U_j,U_l.$
Then for each $t>0$
$$\lim\limits_{n\rightarrow\infty}P\left\{n^3\left(\frac{3\sqrt{3}}{4}-G_n\right)\leq
t\right\}=1-\exp\left\{-\frac{2t}{9\pi}\right\}.$$ The rate of
convergence is $O\left(n^{-\frac{1}{2}}\right)$.
\end{theorem}

Commenting this result, we note that the area of the triangle
inscribed into the unit circumference has the maximal value $
\frac{3\sqrt{3}}{4} $ when its vertices are the vertices of
regular triangle \cite{yaglom}. The following Lemma
\ref{lemma_B} plays an important role in the proof of
Theorem \ref{theorem_B}.

\begin{lemma}
\label{lemma_B} Let $U_1, U_2, U_3$ be independent and uniformly
distributed points on the unit circumference $\mathbb{S}.$ Then
$$\lim\limits_{s\rightarrow +0} s^{-1}P\left\{area(U_1,U_2,U_3)\geq
\frac{3\sqrt{3}}{4}-s\right\}=\frac{4}{3\pi}.$$
\end{lemma}

In  this paper we consider the limit behavior of more general
$U$-max - statistics of this type related to $m$-polygons, $m \geq 3.$
In the sequel $C_1, C_2,\dots$ denote positive constants depending only on $m.$

\section{Inscribed polygon}
\subsection{Perimeter of inscribed polygon}
\renewcommand{\thetheorem}{\arabic{theorem}}
\renewcommand{\thelemma}{\thetheorem.\arabic{lemma}}
\setcounter{theorem}{0}
\setcounter{lemma}{0}

The result of this section is  the genera\-li\-za\-tion of
Theorem \ref{theorem_A} for inscribed triangles to the case of
convex $m$-polygons, $m \geq 3.$ We underline that the proof of this theorem
in \cite{lao} turned out to be inapplicable for random perimeters and areas of
$m$-polygons with $m>3.$  Therefore we had to use some new
ideas.

\begin{theorem}
\label{theorem_1} Let $U_1, U_2, \ldots$ be independent and
uniformly distributed points on the unit circumference
$\mathbb{S}$, and let
\begin{equation}
\label{perim}
 {\mathcal P}_n^m=\max\limits_{1\leq
i_1<\ldots<i_m\leq n} peri(U_{i_1},\ldots,U_{i_m})
\end{equation}
be the maximal perimeter among the perimeters $ {\mathcal P}_m$ of all convex $m$-polygons, generated by $m$ points
from $U_1,\ldots,U_n, \ m \ge 3$. Then for each $t>0$ we have
$$\lim\limits_{n\rightarrow\infty}P\left\{n^{\frac{2m}{m-1}}
\left(2m\sin\frac{\pi}{m} -{\mathcal P}_n^m\right)\leq
t\right\}=1-\exp\left\{-\frac{t^{\frac{m-1}{2}}}{K_{1m}}\right\},$$
where
$K_{1m}=m^{\frac{3}{2}}\left(\pi\sin\frac{\pi}{m}\right)^{\frac{m-1}{2}}
\Gamma\left(\frac{m+1}{2}\right).$ The rate of convergence is
$O\left(n^{- \frac{1}{m-1}}\right)$.
\end{theorem}
Next Lemma is crucial in the proof of Theorem 1.

{\bf Lemma 2.1.}{\it  Let $U_1,\dots, U_m$ be  $m$ independent and uniformly distributed
points  on the unit circumference.
Consider the convex inscribed $m$-polygon with such vertices and with  the perimeter
${\mathcal P}^m=\emph{peri}(U_1, \ldots, U_m).$ We have the
following limit relation:
\begin{equation}
\label{lem1}
\lim\limits_{s\rightarrow +0} s^{-\frac{m-1}{2}}P\left\{{\mathcal P}^m \geq 2m\sin\frac{\pi}{m} -
s\right\}=\frac{\Gamma(m+1)}{K_{1m}},
\end{equation}
where the constant $K_{1m}$ is from
Theorem \ref{theorem_1}.}

\textit{Proof of Lemma 2.1}. For $m=3$ Lemma 2.1
coincides with the result of Lao and Mayer
~\cite{lao}. The perimeter ${\mathcal P}^m  $ is maximal ~\cite{yaglom} for the regular
$m$-polygon (its side is $2\sin\frac{\pi}{m}$) and  this maximal value equals
$2m\sin\frac{\pi}{m}.$ Consider  for $i= 1, \ldots, m-1,
\,  m\ge 3$, the central angle $\beta_i=\angle U_1OU_{i+1}$. By
rotational symmetry, these angles are independent and uniformly
distributed on the interval $[0,2\pi]$. Let show that  the
points $U_i$ can be taken in order of increasing angles, so that
our $m$-polygon corresponds to the following figure.

\begin{center}
\includegraphics[width=6cm]{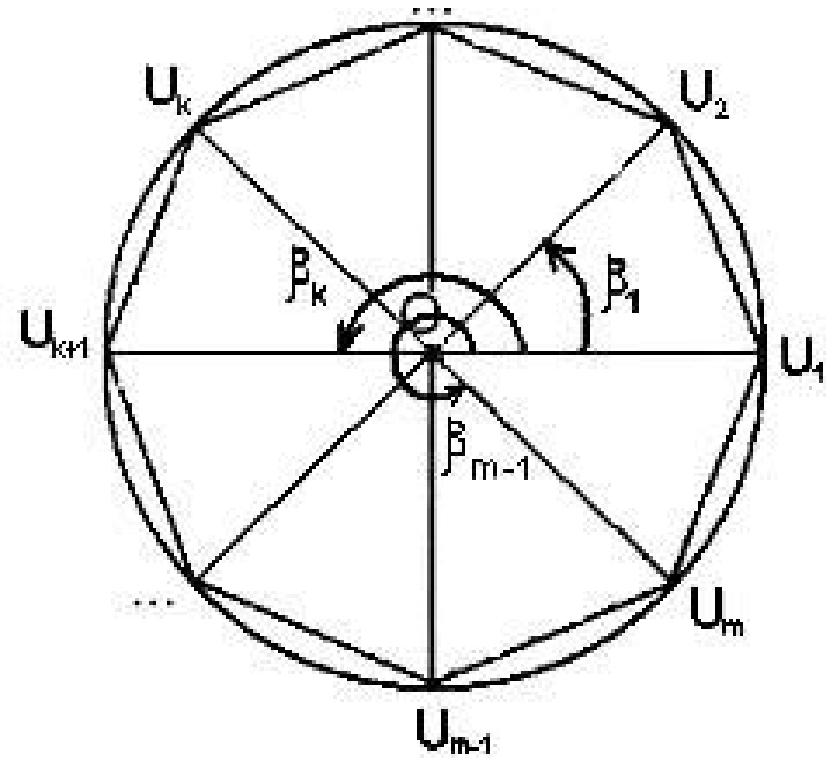}\\
Fig. 1\\
\end{center}

Consider the permutation of angles with increasing indices
$\{\beta_1\leq\beta_2\leq\ldots\leq\beta_{m-1}\}$. Let prove  for any $0 < z < 2m\sin\frac{\pi}{m}$ the equality of probabilities
$$
P\{ {\mathcal P}^m \geq z,
\beta_{i_1}\leq\beta_{i_2}\leq\ldots\leq\beta_{i_{m-1}}\}  = P\{ {\mathcal P}^m \geq z,
\beta_1\leq\beta_2\leq\ldots\leq\beta_{m-1} \},
$$
where $ \{i_1, \ldots, i_{m}\}$ is any permutation of indices $\{1,\ldots,m\}$ and $\beta_{i_s}$ is the central angle defined by the points
$U_{i_{1}} $ and $U_{i_{s +1}}, s = 1,\dots, m-1.$

Let $\mathbb {C}_{m-1}$ be the $(m-1)$-dimensional cube with the edge $2\pi.$ Note that
the sides of $m$-polygon are calculated by the law of cosines:
$$
\label{side_th_cos}
 |U_l U_{l+1}|= 2  \sin \frac{\beta_l -\beta_{l-1}}{2}, 1\leq l \leq m-1; \, |U_m U_1|= 2\sin \frac{\beta_{m-1}}{2}.
$$
We have
\begin{multline*}
P\{ {\mathcal P}^m \geq z,
\beta_{i_1}\leq\beta_{i_2}\leq\ldots\leq\beta_{i_{m-1}}\} = \\ =
P\{\sum\limits_{k=1}^{m-1}|U_{i_k}U_{i_{k+1}}|+|U_{i_{m}}U_{i_{1}}|\geq
z, \beta_{i_1}\leq\beta_{i_2}\leq\ldots\leq\beta_{i_{m-1}} \} = \\
=\frac{1}{(2\pi)^{m-1}} \int_{\mathbb{C}_{m-1}}  {\bf  1}\{
\sum\limits_{k=1}^{m-1}|U_{i_k}U_{i_{k+1}}|+|U_{i_{m}}U_{i_{1}}|
\geq z,  \beta_{i_1}\leq\beta_{i_2}\leq\ldots\leq\beta_{i_{m-1}} \} d\beta_{i_1}\ldots d\beta_{i_{m-1}}.
\end{multline*}
After the change of variables $\beta_{i_1} =
\beta_1,\ldots, \beta_{i_{m-1}} = \beta_{m-1}$   the last integral becomes
\begin{multline*}
 \frac{1}{(2\pi)^{m-1}} \int_{\mathbb{C}_{m-1}}  {\bf  1}\{  |U_{1}U_{2}|+\ldots +|U_{m}U_{1}|  \geq z,
 \beta_{1}\leq\ldots\leq\beta_{m-1} \} d\beta_{1}\ldots d\beta_{m-1} =\\ =
  P\{{\mathcal P}^m \geq z, \beta_{1}\leq\beta_{2}\leq\ldots\leq\beta_{m-1}\},
\end{multline*}
as required.  \hfill$\square$\medskip

Our arguments imply that
\begin{equation}
\label{ineq}
P\{{\mathcal P}^m \geq 2m\sin\frac{\pi}{m}-s\} = (m-1)!P\{{\mathcal P}^m \geq 2m\sin\frac{\pi}{m}-s, \  \beta_{1}\leq\beta_{2}\leq\ldots\leq\beta_{m-1}  \}.
\end{equation}
The condition $\beta_{1}\leq\beta_{2}\leq\ldots\leq\beta_{m-1}$ provides that the points
$U_1,\dots, U_m$ stand in increasing order just as in Fig.1.

Now let show that for small $s>0 $ and under the condition
\begin{equation}
\label{crucial}
{\mathcal P}^m \geq 2m\sin\frac{\pi}{m}-s
\end{equation}
all central angles of our $m$-polygon
differ very little from $\frac{2\pi}{m}.$ Take the smallest  and the largest side of our polygon with opposite central angles
$\frac{2\pi}{m} - 2\psi$ and  $ \frac{2\pi}{m} +2\varphi, \  \varphi > \psi > 0.$ One can displace these sides so that they are adjacent. Their lengths are equal to $2\sin( \frac{\pi}{m} +\varphi)$ and $2\sin( \frac{\pi}{m} - \psi).$ Now let make  the so-called "symmetrization"  \cite{yaglom}  replacing
our two sides by two equal sides having  the common vertex on the circle in the center of the arc subtending the sum of the angles.
Then the length of new  two sides will increase \cite[Probl. 55b]{yaglom} and  is equal to  $4\sin ( \frac{\pi}{m} +\varphi- \psi). $ The increment of the perimeter is just
$$
\Delta = 4 \sin \left( \frac{\pi}{m} +\frac{\varphi- \psi}{2}\right)\left(1-\cos\frac{\varphi+ \psi}{2}\right) , \, \varphi - \psi > 0.
$$
Suppose the initial perimeter of our $m$-polygon was  $2m\sin\frac{\pi}{m}- \sigma$ with $\sigma \leq s.$  Consequently
$$
2m\sin\frac{\pi}{m}- \sigma + \Delta  \leq 2m\sin\frac{\pi}{m},
$$
whence it follows that $ \Delta  \le s $ and therefore
\begin{equation}
\label{cos}
1 - \cos\frac{\varphi+ \psi}{2} \leq  \frac{s}{ 4 \sin \left( \frac{\pi}{m} +\frac{\varphi- \psi}{2} \right) }.
\end{equation}

Clearly, the sum of the largest and smallest central angles $\frac{4\pi}{m} +2\varphi-2\psi$ is smaller than $2\pi$. Hence
$0< \frac{\pi}{m} + \frac{\varphi - \psi}{2} < \frac{\pi}{2}$  and we have from (\ref{cos}) the inequality
$$
1 - \cos\frac{\varphi+ \psi}{2} \leq \frac{ s}{ 4 \sin  \frac{\pi}{m}}  =  C_1 s,
$$
which implies that for small $s> 0$
$$
\varphi+ \psi \le 2\arccos(1- C_1s) < C_2 \sqrt{s}.
$$
Consequently all random central angles differ from the angle $\frac{2\pi}{m}$ by no more than $O(\sqrt{s}.)$

Let estimate the deviation of the angles $\beta_k$ from $\beta_{k-1}, k=2,
\ldots, m-1,$ under  the condition (\ref{crucial}).
We introduce the auxiliary random angles $\alpha_0=0, \alpha_1, \ldots, \alpha_{m-1} , \alpha_m =0$ such that
\begin{equation}
\label{beta_k}
\beta_k=\frac{2\pi k}{m}+\alpha_k,\ k=1, \ldots,m.
\end{equation}
The random angles $ \alpha_1, \ldots, \alpha_{m-1} $ are independent and each $\alpha_k$ is uniformly distributed on $\left[-\frac{2\pi
k}{m},2\pi-\frac{2\pi k}{m}\right].$  In terms of $\alpha_k$  we have
$$
{\mathcal P}^m=2\sum\limits_{k=1}^{m}\sin\left(\frac{\pi}{m}+\frac{\alpha_k-\alpha_{k-1}}{2}\right),
$$
 and the inequality (\ref{crucial}) takes the form
\begin{equation}
\label{sum_sin}
2\sum\limits_{k=1}^{m}\sin\left(\frac{\pi}{m}
+\frac{\alpha_k-\alpha_{k-1}}{2}\right)\geq
2m\sin\frac{\pi}{m}-s.
\end{equation}

The argument given above shows that any central angle differs from the expected value $\frac{2\pi}{m}$ by no more than $O(\sqrt{s}.)$ Hence we have under (\ref{crucial})
\begin{equation}
\label{differences}
\max_{1\le k \le m} |\alpha_k - \alpha_{k-1}| \le  C_3\sqrt{s}.
\end{equation}
and consequently
\begin{equation}
\label{alphas}
\max_{1\le k \le m-1} |\alpha_k | \le  C_4\sqrt{s}.
\end{equation}

Let return to  the formula (\ref{sum_sin}).
As the differences $|\alpha_k-\alpha_{k-1}|$ are small, we can expand in (\ref{sum_sin}) the sine function in Taylor series with the remainder term. Hence  we obtain for some small random angles $\eta_k, k=1,\ldots,m-1,$ that
\begin{multline*}
2\sum\limits_{k=1}^{m}\sin\left(\frac{\pi}{m}+\frac12 (\alpha_k-\alpha_{k-1}) \right)=\\ = 2m\sin\frac{\pi}{m}-
\frac{1}{4}\sin\frac{\pi}{m}\sum\limits_{k=1}^{m}(\alpha_k-\alpha_{k-1})^2+ \frac{1}{24}\sum\limits_{k=1}^{m}\cos\left(\frac{\pi}{m} +\eta_k\right)(\alpha_k-\alpha_{k-1})^3,
\end{multline*}
and therefore (\ref{sum_sin}) is equivalent to
\begin{equation}
\label{long}
\sum\limits_{k=1}^{m}(\alpha_k-\alpha_{k-1})^2 - \frac{1}{6\sin \frac{\pi}{m}}\sum\limits_{k=1}^{m}\cos\left(\frac{\pi}{m} +\eta_k\right)(\alpha_k-\alpha_{k-1})^3 \leq  \frac{4s}{ \sin \frac{\pi}{m} }.
\end{equation}

Clearly
$$
|\sum\limits_{k=1}^{m}\cos\left(\frac{\pi}{m} +\eta_k\right)(\alpha_k-\alpha_{k-1})^3| \leq \max_{1\leq k \leq m} |\alpha_k -\alpha_{k-1}|  \sum\limits_{k=1}^{m}(\alpha_k-\alpha_{k-1})^2,
$$
hence,  using  (\ref{long}) and (\ref{differences}), we may write the inequalities
\begin{multline*}
P\left\{{\mathcal P}^m \geq 2m\sin\frac{\pi}{m}-s, \  \beta_{1}\leq\beta_{2}\leq\ldots\leq\beta_{m-1}  \right\} \le\\
\le P\left\{ \sum\limits_{k=1}^{m}(\alpha_k-\alpha_{k-1})^2  \le \frac{24s}{ 6  \sin \frac{\pi}{m} - \max_{1\le k \le m} |\alpha_k-\alpha_{k-1}|} \right \}\le\\  \le P\left\{\sum\limits_{k=1}^{m}(\alpha_k-\alpha_{k-1})^2 \leq \frac{4s_1}{\sin \frac{\pi}{m}}\right\},
\end{multline*}
where $  s_1 \to s$ as $s \to 0.$
Quite analogously
$$
P\left\{{\mathcal P}^m \geq 2m\sin\frac{\pi}{m}-s, \  \beta_{1}\leq\beta_{2}\leq\ldots\leq\beta_{m-1}  \right\} \ge
P\left\{\sum\limits_{k=1}^{m}(\alpha_k-\alpha_{k-1})^2 \leq \frac{4s_2}{\sin \frac{\pi}{m}}\right\},
$$
where $ s_2 \to s$ as $s \to 0.$

Now let introduce the quadratic form $Q(\alpha) = Q(\alpha_1,\dots,\alpha_{m-1})$ by
$$
Q(\alpha) = \frac12 \sum\limits_{k=1}^{m}(\alpha_k- \alpha_{k-1})^2 =  \sum\limits_{k=1}^{m-1}\alpha_k^2-  \sum\limits_{k=1}^{m-1}\alpha_k\alpha_{k-1}.
$$
We have by (\ref{ineq})  for small $s_1 > s_2 , s_1 \to s , s_2 \to s$  as $s \to 0:$
\begin{equation}
\label{main}
 P\left\{Q\left(\alpha\right) \leq
\frac{2s_2}{ \sin\frac{\pi}{m} }\right\} \leq  \frac{1}{ (m-1)!} P\left\{{\mathcal P}^m \geq 2m\sin\frac{\pi}{m}-s\right\} \leq  P\left\{Q\left(\alpha\right) \leq
\frac{2s_1}{ \sin\frac{\pi}{m} }\right\}.
\end{equation}

Now we proceed  to the calculation of the probability in the right-hand side of (\ref{main}).
The quadratic form  $Q(\alpha)$  has the following matrix of size
$(m-1)\times(m-1)$:
$$B=\left(\begin{array}{ccccccc}
    1  & -1/2 &   0  & 0    & 0    & \ldots & 0 \\
  -1/2 &   1  & -1/2 & 0    & 0    & \ldots & 0 \\
    0  & -1/2 &   1  & -1/2 & 0    & \ldots & 0 \\
    0  &   0  & -1/2 & 1    & -1/2 & \ldots & 0 \\
    \vdots & \vdots & \vdots & \ddots & \ddots & \ddots &  \vdots\\
    0 & \ldots & \ldots  &   \ldots &   -1/2 &      1 &   -1/2 \\
    0 & \ldots & \ldots &     \ldots &     0  &   -1/2 & 1
  \end{array}\right).
$$
This matrix is symmetric tridiagonal. The spectrum of such matrices is
known, see \cite[p.137]{spectrum} or \cite{taiwan}, so that all eigenvalues of the matrix $B$ are\\
$$\lambda_k=1-\cos \frac{\pi k}{m}, \quad k=1,\ldots,m-1.$$

Then using the appropriate orthogonal transformation, we can replace
the quadratic form
$\sum\limits_{k=1}^{m-1}\alpha_k^2-\sum\limits_{k=1}^{m-1}\alpha_k\alpha_{k-1}$
by the  quadratic form
\begin{equation} \label{orthogonal_transformation}
\sum\limits_{k=1}^{m-1}\left(1-\cos \frac{\pi
k}{m}\right)Y_k^2\\
\end{equation}
in new variables $Y_k.$ Now set for brevity
$$ l_k^2 = 2s_1/ \sin  \frac{\pi}{m} \left(1-\cos \frac{\pi k}{m}\right), k = 1, \ldots, m-1 .$$

As the  Jacobian of the orthogonal  transformation is 1, we get by (\ref{main}) that
\begin{multline*}
P\left\{{\mathcal P}^m \geq 2m\sin\frac{\pi}{m} -
s\right\} \le \frac{(m-1)!}{(2\pi)^{m-1}}\int_{\left\{\sum\limits_{k=1}^{m-1}x_k^2-\sum\limits_{k=1}^{m-1}x_k x_{k-1}
\leq\frac{2s_1}{\sin\frac{\pi}{m}}\right\}} \ dx_1\ldots
dx_{m-1}= \\ =
\frac{(m-1)!}{(2\pi)^{m-1}}\int_{\left\{\sum\limits_{k=1}^{m-1}\left(1-\cos\frac{\pi
k}{m}\right)y_k^2\leqslant
\frac{2s_1}{\sin\frac{\pi}{m}}\right\}} \ dy_1\ldots
dy_{m-1} = \\ =\frac{(m-1)!}{(2\pi)^{m-1}}
\int_{\left\{\sum\limits_{k=1}^{m-1}\left(\frac{y_k}{l_k }\right)^2 \leqslant 1\right\}} \  dy_1\ldots
dy_{m-1}= \frac{(m-1)!}{(2\pi)^{m-1}}b_{m-1} \Pi_{m-1}(s_1).
\end{multline*}
Here $b_{m-1}$ is a  volume of $(m-1)$-dimensional ball of unit
radius, and $\Pi_{m-1}(s_1)$ is a product of semi-axes of the small
ellipsoid
$$
\left\{(y_1,\dots,y_{m-1}):\sum\limits_{k=1}^{m-1}\frac{y_k^2}{l_k^2}\leqslant 1\right\} ,
$$
which is equal to
$$\Pi_{m-1}(s_1)=\prod\limits_{k=1}^{m-1} l_k = \prod\limits_{k=1}^{m-1}\sqrt{\frac{2s_1}{\sin\frac{\pi}{m}\left(1-\cos\left(\frac{\pi
}{m}\right)\right)}}=\left(\frac{s_1}{\sin\frac{\pi}{m}}\right)^{\frac{m-1}{2}}\prod\limits_{k=1}^{m-1}
\left(\sin\left(\frac{\pi k}{2m}\right)\right)^{-1}.$$

\noindent To simplify this expression we use the identity from \cite[formula (6.1.2.3)]{Prud}:
$$
\prod\limits_{k=1}^{m-1}\sin\left(\frac{\pi
k}{2m}\right)=\frac{\sqrt{m}}{2^{m-1}}.
$$
Hence the product of semi-axes of the ellipsoid  is equal to
$$\Pi_{m-1}(s_1)=\left(\frac{s_1}{\sin\frac{\pi}{m}}\right)^{\frac{m-1}{2}}\frac{2^{m-1}}{\sqrt{m}}=
\frac{2^{m-1}s_1^{\frac{m-1}{2}}}{\sqrt{m}\left(\sin\frac{\pi}{m}\right)^{\frac{m-1}{2}}}.$$
The volume of $(m-1)$-dimensional ball of the unit radius is
well-known  and equals
\begin{equation}
\label{volume_b_m}
b_{m-1}=\frac{\left(\sqrt{\pi}\right)^{m-1}}{\Gamma\left(\frac{m+1}{2}\right)}.
\end{equation}

Therefore the volume of the $(m-1)$-dimensional ellipsoid is equal to
$$V_{m-1}(s_1) =b_{m-1}\Pi_{m-1}(s_1)=\frac{\left(\sqrt{\pi}\right)^{m-1}}{\Gamma\left(\frac{m+1}{2}\right)}
\frac{2^{m-1}s_1^{\frac{m-1}{2}}}{\sqrt{m}\left(\sin\frac{\pi}{m}\right)^{\frac{m-1}{2}}} ,$$
from where we finally obtain:
$$P\left\{{\mathcal P}^m\geq
2m\sin\frac{\pi}{m}-s\right\} \le  \frac{(m-1)!}{(2\pi)^{m-1}}V_{m-1}(s_1)=\frac{(m-1)!}{(2\pi)^{m-1}}
\frac{\left(\sqrt{\pi}\right)^{m-1}}{\Gamma\left(\frac{m+1}{2}\right)}\frac{2^{m-1}
s_1^{\frac{m-1}{2}}}{\sqrt{m}\left(\sin\frac{\pi}{m}\right)^{\frac{m-1}{2}}}=$$
$$=\frac{s_1^{\frac{m-1}{2}}}{\sqrt{m} \left(\pi \sin\frac{\pi}{m}\right)^{\frac{m-1}{2}}}
\frac{\Gamma(m)}{\Gamma\left(\frac{m+1}{2}\right)}=s_1^{\frac{m-1}{2}}\frac{\Gamma(m+1)}{K_{1m}}.$$
In the same manner we get from (\ref{main}) the opposite inequality
$$P\left\{{\mathcal P}^m\geq
2m\sin\frac{\pi}{m}-s\right\} \ge  s_2^{\frac{m-1}{2}}\frac{\Gamma(m+1)}{K_{1m}}.$$
Lemma 2.1 immediately follows from these two inequalities.  \hfill$\square$\medskip

\textit{Proof of Theorem \ref{theorem_1}.}
Consider the transformation $$z_n(t)=2m\sin\frac{\pi}{m}-s = 2m\sin\frac{\pi}{m}-tn^{-\frac{2m}{m-1}},$$
and denote
$$\lambda_{n,z_n(t)}=\binom{n}{m}P\{{\mathcal P}_n^m>z_n(t)\}=\frac{n!}{m!(n-m)!}P\{{\mathcal
P}_n^m>z_n(t)\},$$
where $ {\mathcal P}_n^m$ is the maximal perimeter from (\ref{perim}).

As $s=t n^{-\frac{2m}{m-1}},$ we note that  $n^m s^{\frac{m-1}{2}}=t^{\frac{m-1}{2}}.$
Since for fixed $m$ and  large $n$ it holds  $\frac{n!}{(n-m)!}\sim
n^m$, then we have from Lemma 2.1 that
\begin{multline*}
\lim\limits_{n\rightarrow\infty}
\lambda_{n,z_n(t)}=\frac{1}{m!}\lim\limits_{n\rightarrow\infty}n^ms^{\frac{m-1}{2}}\cdot
s^{-\frac{m-1}{2}}P\{{\mathcal P}_n^m>2m\sin\frac{\pi}{m}-s\} =\\ =\frac{1}{m!}t^{\frac{m-1}{2}}\lim\limits_{n\rightarrow\infty}\left(tn^{-\frac{2m}{m-1}}\right)^{-\frac{m-1}{2}}P\{{\mathcal
P}_n^m>2m\sin\frac{\pi}{m}-tn^{-\frac{2m}{m-1}}\}= \\ = \frac{1}{m!}t^{\frac{m-1}{2}}\frac{\Gamma(m+1)}{K_{1m}}=\frac{t^{\frac{m-1}{2}}}{K_{1m}}:=\lambda_t>0.
\end{multline*}
Thus, the condition (\ref{corollary_silv_1}) from Silverman-Brown Theorem  holds true.

Now for any $ r \in \{1,\ldots,m-1\} $ denote by ${\mathcal P}^{m;i_1,\ldots, i_m} $ the perimeter of the convex inscribed $m$-polygon based on
the random points $U_{i_1}, \ldots, U_{i_m}, i_1< \ldots <i_m,$ on the unit circumference.  We must verify the condition
(\ref{corollary_silv_4}):
$$
\lim\limits_{n\rightarrow\infty}n^{2m-r}P\{{\mathcal P}^{m;1,\ldots,m}>z_n(t),\, {\mathcal P}^{m; m - r +1,\ldots,2m-r}>z_n(t)\}=0
$$
for each $r\in\{1,\ldots,m-1\}$.

{\bf Lemma 2.2} (verification of the condition (\ref{corollary_silv_4})).
{\it  For each $r\in\{1,\ldots,m-1\}$ it holds}
$$\lim\limits_{n\rightarrow\infty}n^{2m-r}P\{{\mathcal
P}^{m;1,\ldots,m}>z_n(t),{\mathcal P}^{m;
m - r+1,\ldots,2m-r}>z_n(t)\}=0.
$$
\noindent \textit{Proof of Lemma 2.2}.
Using same arguments as in the proof of Lemma 2.1, we may assume that $\beta_1<\beta_2<...<\beta_{m-1},$ so that the points $U_j$ follow
one after one in increasing order. Then the condition ${\mathcal P}^{m;1,\ldots,m}>z_n(t)$ implies by (\ref{alphas})
 that for some constant $C_5>0$
$$
|\alpha_1| < C_5\sqrt{s},  ...  , \  |\alpha_{m-1}| < C_5 \sqrt{s} .
$$
The second condition ${\mathcal P}^{m;m - r+1,\ldots,2m-r}>z_n(t)$  implies analogously that
$$
|\alpha_{m-r+1}| < C_5\sqrt{s}, ...  ,\  |\alpha_{2m-r-1} | < C_5\sqrt{s} .
$$
The intersection of these two events is the event
$$
\bigcap_{j=1}^{2m-r-1} \{ \alpha_j < C_5\sqrt{s} \}.
$$
All angles $\alpha_j$ are independent and have the uniform distribution on the intervals \\ $\left(-\frac{2\pi
j}{m},2\pi-\frac{2\pi j}{m}\right)$  of length $2\pi.$
Then the expression of interest can be estimated for each $r=1,...,m-1, $ as follows
\begin{multline*}
n^{2m-r}P\{{\mathcal P}^{m;1,\ldots,m}>z_n(t),{\mathcal P}^{m; m - r+1,\ldots,2m-r}>z_n(t)\} \leq \\
 \\ \leq n^{2m-r} \left( C_5 \sqrt{s}/2\pi \right)^{2m - r -1}= O( n^{2m-r - \frac{(2m - r -1)m}{m-1}}) = O(n^{- \frac{m-r}{m-1}}) = o(1).
 \end{multline*}
 \hfill$\square$

By Remark \ref{remark_3} it follows that the rate of
convergence in the worst case  $r = m-1$ is
$$O\left(n^{m+1}p_{n,z_n(t)}\tau_{n,z_n(t)}(m-1)\right)=O\left(n^{m+1}n^{-\frac{m^2}{m-1}}\right)=O\left(n^{-\frac{1}{m-1} }\right).$$

We see that this result coincides with the result of Lao and Mayer in Theorem \ref{theorem_A} for $m=3$ but deteriorates when $m$ grows.

\textit{Continuation of the proof of Theorem \ref{theorem_1}.}
Next  we apply the relation (\ref{corollary_silv_3}), see also \cite{lao}, and obtain:
$$\lim\limits_{n\rightarrow\infty}P\left\{{\mathcal P}_n^m<2m\sin\frac{\pi}{m} -tn^{-\frac{2m}{m-1}}\right\}=
\exp\left\{-\frac{t^{\frac{m-1}{2}}}{K_{1m}}\right\}.$$

Consequently  for each $t>0$ we have:
$$\lim\limits_{n\rightarrow\infty}P\left\{n^{\frac{2m}{m-1}}(2m\sin\frac{\pi}{m}
-{\mathcal P}_n^m)\leq
t\right\}=1-\exp\left\{-\frac{t^{\frac{m-1}{2}}}{K_{1m}}\right\},$$
where
$$K_{1m}=m^{\frac{3}{2}}\left(\pi\sin\frac{\pi}{m}\right)^{\frac{m-1}{2}}
\Gamma\left(\frac{m+1}{2}\right).$$
\hfill$\square$

The limit distribution is the Weibull distribution which is
expected for the distribution of a maximum. For $m=3$ we get the result from \cite{lao},
i.e. Theorem \ref{theorem_A} from Introduction.

\subsection{Area of inscribed polygon}

\begin{theorem}
\label{theorem_2} Let $U_1, U_2, \ldots$ be independent and
uniformly distributed points on the unit circumfe\-ren\-ce  $\mathbb{S},\ {\mathcal A}_n^m=\max\limits_{1\leq
i_1<\ldots<i_m\leq n} {\mathcal A}_m(U_{i_1},\ldots,U_{i_m})$ is a
maximal area  among all areas of convex $m$-polygons with the vertices
$U_{i_1},\ldots,U_{i_m}$. \ Then for each $t>0$ it is true that
$$\lim\limits_{n\rightarrow\infty}P\left\{n^{\frac{2m}{m-1}}
\left(\frac{m}{2}\sin\frac{2\pi}{m} -{\mathcal A}_n^m\right)\leq
t\right\}=1-\exp\left\{-\frac{t^{\frac{m-1}{2}}}{K_{2m}}\right\},$$
where
\begin{equation}
\label{const_2}
K_{2m}=m^{\frac{3}{2}}\left(\pi\sin\frac{2\pi}{m}\right)^
{\frac{m-1}{2}}\Gamma\left(\frac{m+1}{2}\right).
\end{equation}
The rate of
convergence is $O\left(n^{- \frac{1}{m-1}}\right)$.
\end{theorem}

\begin{rem}
\label{remark_5} Curiously, the limit constants $K_{1m} $ and
$K_{2m}$ are very similar, but yet differ by the argument of the sine function. At
the end of the paper we will discuss the asymptotic behavior of
them and of similar constants from limit theorems for the
metric characteristics of circumscribed polygons.
\end{rem}

{\bf Lemma 2.3}. {\it Let $m$ independent and uniformly distributed
points $U_1, \ldots, U_m$ be chosen on the unit circumference.
Consider the convex $m$-polygon with such vertices having the area ${\mathcal
A}^m=area(U_1, \ldots, U_m).$ Then the following limit relation holds:}
$$\lim\limits_{s\rightarrow +0} s^{-\frac{m-1}{2}}P\left\{{\mathcal
A}^m \geq
\frac{m}{2}\sin\frac{2\pi}{m}-s\right\}=\frac{\Gamma(m+1)}{K_{2m}},$$
where $K_{2m}$ is given by (\ref{const_2}).

\textit{Proof of Lemma 2.3}. The area of $m$-polygon
inscribed into the unit circle has the maximal value for the
regular $m$-polygon \cite{yaglom} and equals $\max {\mathcal
A}^m=\frac{m}{2}\sin\frac{2\pi}{m}$. As in the proof of Lemma
2.1, consider the angles $\beta_i$.

 The area of
inscribed $m$-polygon is the sum of the areas $S_i$ of triangles,
formed by the triples of points  $U_i, U_{i+1},O,$
where $O$ is the center of the circle, see Fig.2.

\begin{center}
\includegraphics[width=7cm]{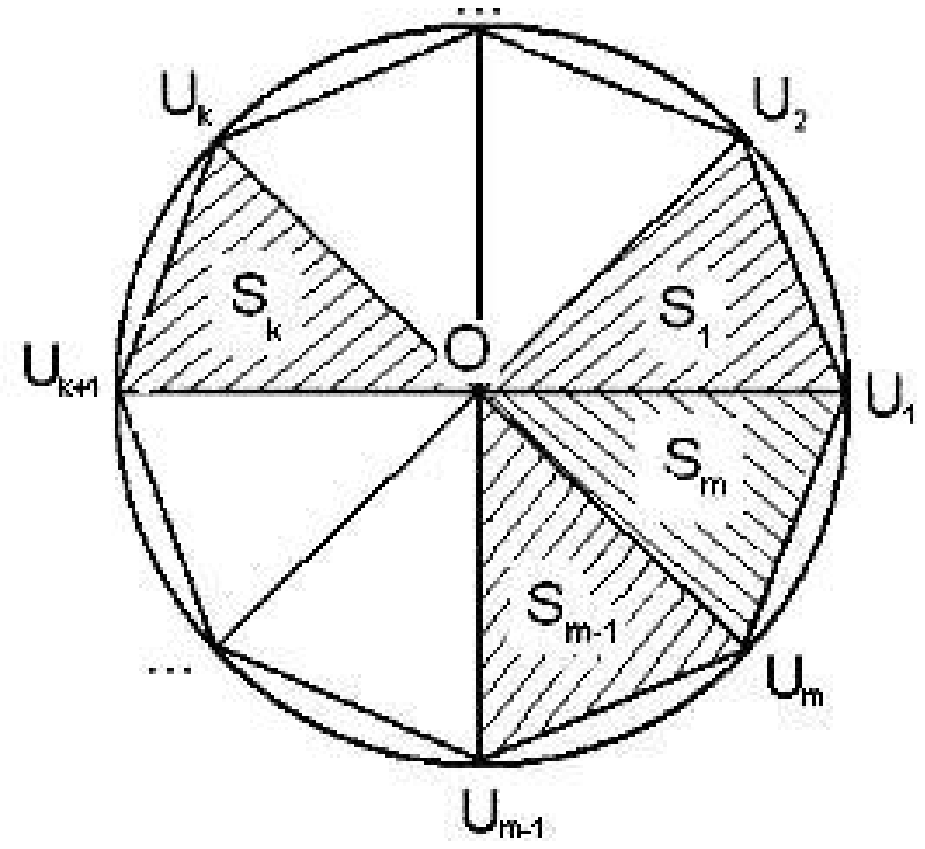}\\
Fig. 2\\
\end{center}

Given that the angle between $U_k$ and $U_{k+1}$ equals
$\beta_k-\beta_{k-1},$  the area of the triangle $\triangle
U_kOU_{k+1}$ is $\frac{1}{2}\sin(\beta_k-\beta_{k-1})$. Hence
$${\mathcal A}^m=\sum\limits_{k=1}^{m}S_k=\frac{1}{2}\sum\limits_{k=1}^{m}\sin(\beta_k-\beta_{k-1}).$$
We have, assuming that $\beta_0=0, \beta_m=2\pi,$
\begin{multline*}
P\left\{{\mathcal A}^m\geq\frac{m}{2}\sin\frac{2\pi}{m}-s\right\}=\\
=(m-1)!P\left\{\frac{1}{2}\sum\limits_{k=1}^{m}\sin(\beta_k-\beta_{k-1})\geq\frac{m}{2}\sin\frac{2\pi}{m}-s,
\beta_1\leq\beta_2\leq\ldots\beta_{m-1}\right\}.
\end{multline*}

Next we pass by (\ref{beta_k}) from random angles $\beta_k$ to random angles $\alpha_k,$
which are independent and uniformly distributed on $ [ -\frac{2\pi k}{m}, 2\pi -  \frac{2\pi k}{m}]$.

Same arguments as in the proof of Lemma 2.1 show that for small $s$ and under the condition
$$
{\mathcal A}^m\geq\frac{m}{2}\sin\frac{2\pi}{m}-s
$$
all central angles of the inscribed polygon differ from $\frac{2\pi}{m}$ at most by $O(\sqrt{s})$ and
 the  inequalities    (\ref{differences}) and (\ref{alphas})   hold.  Expanding the sine function in the expression of area
 for small $\alpha_k, k = 1,\ldots, m-1$, we get in the same way as above the following
 inequalities for small $s:$
\begin{multline*}
P\left\{\sum\limits_{k=1}^{m-1}\alpha_k^2-\sum\limits_{k=1}^{m-1}\alpha_k\alpha_{k-1}\leq\frac{2s_4}{\sin\frac{2\pi}{m}}\right\} \le\\
\le P\left\{\frac{1}{2}\sum\limits_{k=1}^{m}\sin(\beta_k-\beta_{k-1})\geq\frac{m}{2}\sin\frac{2\pi}{m}-s,
\beta_1\leq\beta_2\leq\ldots\beta_{m-1}\right\} \le\\
\le P\left\{\sum\limits_{k=1}^{m-1}\alpha_k^2-\sum\limits_{k=1}^{m-1}\alpha_k\alpha_{k-1}\leq\frac{2s_3}{\sin\frac{2\pi}{m}}\right\},
\end{multline*}
where $s_3>s_4$, and $s_3\to s, s_4 \to s$ as $s \to 0.$

The probabilities in the left and in the right have been calculated above
(with another constant in the right-hand side), so we get analogously
$$
s_4^{\frac{m-1}{2}}\frac{\Gamma(m+1)}{K_{2m}} \le
P\left\{{\mathcal
A}^m\geq\frac{m}{2}\sin\frac{2\pi}{m}-s\right\} \le s_3^{\frac{m-1}{2}}\frac{\Gamma(m+1)}{K_{2m}},
$$ where
$$K_{2m} = m^{\frac{3}{2}}\left(\pi\sin\frac{2\pi}{m}\right)^
{\frac{m-1}{2}}\Gamma\left(\frac{m+1}{2}\right).
$$
Lemma 2.3 follows from this as $s\to 0.$ \hfill$\square$

 \medskip

\textit{Proof of Theorem \ref{theorem_2}} \, is quite analogous to
the proof of Theorem \ref{theorem_1}, but with the transfor\-mation \
$z_n(t)=\frac{m}{2}\sin\frac{2\pi}{m}-tn^{-\frac{2m}{m-1}}.$ \hfill $\square$

\section{Circumscribed polygon}

We proceed to the perimeter and area of random
{\it circumscribed} polygons. In this case we are of course
interested in {\it minimal} perimeter and area. As far as we
know this problem has never been studied.

\subsection{Perimeter of circumscribed polygon }

\begin{theorem}
\label{theorem_3} Let $U_1,U_2,\ldots$ be independent and
uniformly distributed points on the circum\-fe\-rence
$\mathbb{S}$ of unit radius. For each set of points
$U_{i_1},\ldots,U_{i_m}$ on the circle draw the tangents at
these points, denoting the intersection of these tangents at the points
$U_{i_k}$ and $U_{i_{k+1}}$ by the points $V_{i_k}, k=1,\ldots,m,$
correspondingly. Let $peri_n^m=\min\limits_{1\leq
i_1<\ldots<i_m\leq n} peri_m(V_{i_1},\ldots,V_{i_m})$ be minimal
of perimeters among all circumscribed $m$-polygons with the vertices $V_{i_1},\ldots,V_{i_m}$. Then for each $t>0$ we have
$$\lim\limits_{n\rightarrow\infty}P\left\{n^{\frac{2m}{m-1}}
\left(peri_n^m-2m\tan\frac{\pi}{m}\right)\leq
t\right\} = 1-\exp\left\{-\frac{K_{3m}t^{\frac{m-1}{2}}}{\Gamma(m+1)}\right\},$$
where
$$K_{3m}=m^{\frac{3}{2}}\left(2\pi\left(1+\tan^2\frac{\pi}{m}\right)\tan\frac{\pi}{m}\right)^{\frac{m-1}{2}}
\Gamma\left(\frac{m+1}{2}\right).$$ The rate of convergence is
$O\left(n^{-\frac{1}{m-1}}\right)$.
\end{theorem}

Clearly, the minimal perimeter is achieved at the
regular circumscribed $m$-polygon, whence the centering constant
$2m\tan\frac{\pi}{m}$ appears. The following Lemma 3.1
plays a key role in a proof of Theorem 3.

{\bf Lemma 3.1}.{\it  Let  choose $m$ independent and uniformly
distributed points $U_1, \ldots, U_m$ on the unit circle.
Draw the tangents at these points on the circle and denote the
points of intersection of the tangents
 by $V_k, k=1,\ldots, m,$ correspondingly. Consider the convex
circumscribed $m$-polygon with the vertices $V_1, \ldots, V_m,$ and perimeter $peri_m(V_1,
\ldots, V_m).$  Then the limit relation holds:}
$$\lim\limits_{s\rightarrow +0}
s^{-\frac{m-1}{2}}P\left\{peri_m(V_1,\ldots,V_m)\leq
2m\tan\frac{\pi}{m}+s\right\}=\frac{\Gamma(m+1)}{K_{3m}},$$ where
$$K_{3m} = m^{\frac{3}{2}}\left(2\pi\left(1+\tan^2\frac{\pi}{m}\right)\tan\frac{\pi}{m}\right)^{\frac{m-1}{2}}
\Gamma\left(\frac{m+1}{2}\right).$$

\quad\textit{Proof of Lemma 3.1.}  In the sequel denote for brevity $peri_m(V_1,\ldots,V_m):= peri_m.$
Consider the angles $\beta_k, k=1,\ldots,m-1$, which are uniformly
distributed on $[0,2\pi]$, where
$$\angle U_1OU_k=\beta_{k-1}, k=2,\ldots,m-1, \quad \angle U_mOU_1=2\pi-\beta_{m-1}.$$
From this we get $\angle U_kOU_{k+1}=\beta_k-\beta_{k-1}.$  Next we obtain  by same
arguments as in the proof of Lemma  2.1 that
$$
\label{sum_{m-1}}
P\left\{peri_m \leq 2m \tan \frac{\pi}{m}+s\right\}
=(m-1)!P\left\{peri_m \leq
2m\tan\frac{\pi}{m}+s, \beta_1\leq \ldots \leq \beta_{m-1}\right\}.
$$

Therefore we may assume that the points $U_k$ and therefore
the points $V_k$ go in increasing order.  Denote the segments of the tangents to the circumference
drawn from the point $V_k$ by:
$$|U_kV_k|=|V_kU_{k+1}|=:a_k, k=1,\ldots,m-1, \, |U_mV_m|=|V_mU_1|=:a_m.$$

\begin{center}
\includegraphics[width=7cm]{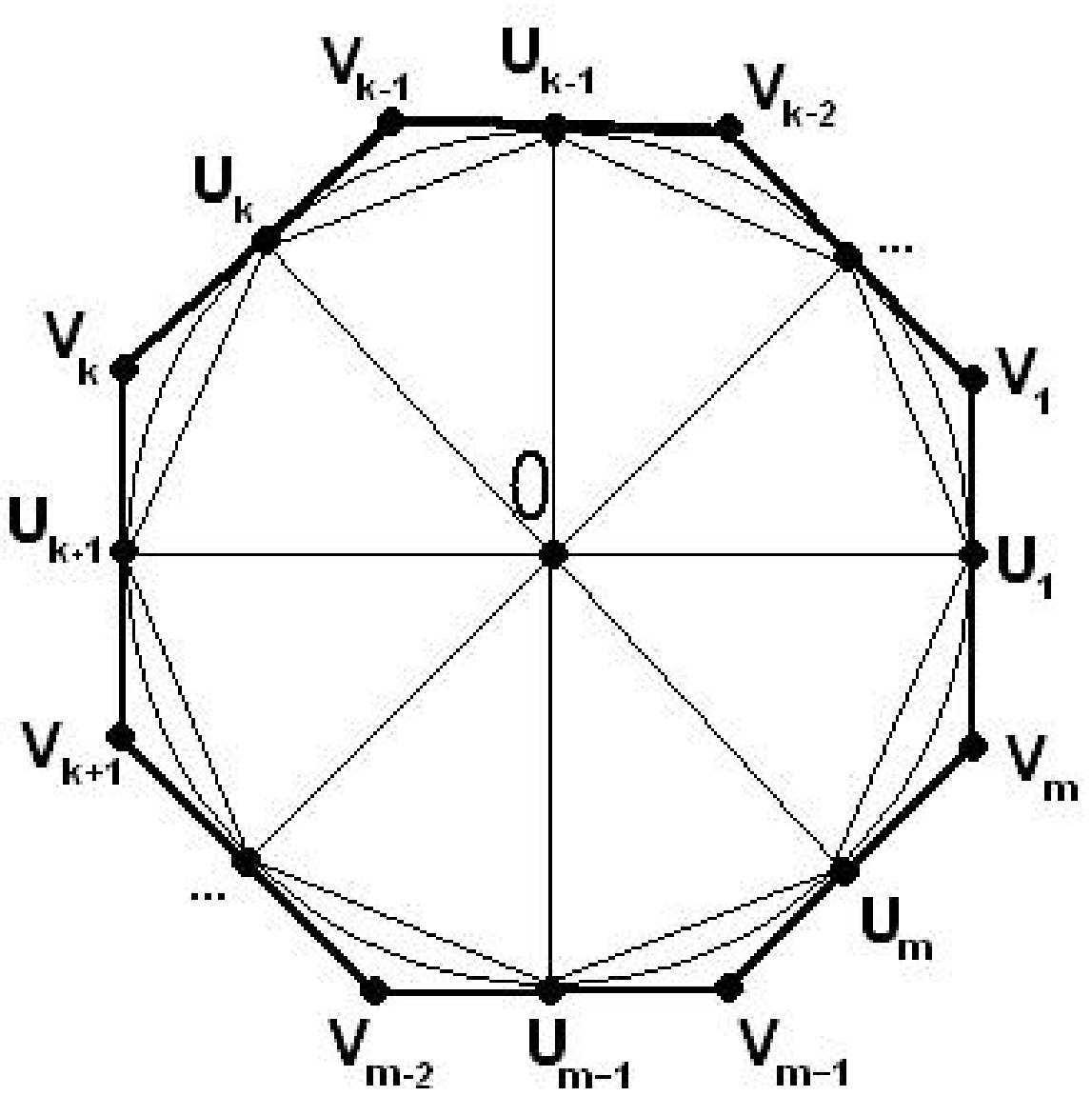}\\
Fig. 3\\
\end{center}

Hence  the perimeter of circumscribed $m$-polygon is given by
$$peri_m = peri_m(V_1, \ldots, V_m) = 2\sum\limits_{k=1}^{m}a_k.$$
Performing some elementary calculations, we obtain
\begin{equation}
\label{a_k}
 a_k=\tan \left(\frac{\beta_k-\beta_{k-1}}{2}\right),\, peri_m=2\sum\limits_{k=1}^{m}\tan\left(\frac{\beta_k-\beta_{k-1}}{2}\right).
\end{equation}
Analogously to the  proof of Lemma 2.1, we pass from the angles $\beta_k$
to the auxiliary angles $\alpha_k, k=1,\ldots,m-1,$ by (\ref{beta_k}):
$$\frac{\beta_k-\beta_{k-1}}{2}=\frac{\pi}{m}+\frac{\alpha_k-\alpha_{k-1}}{2}.$$

As in the proof of Lemma 1.1, one can show that under the condition
$$
 peri_m\leq
2m\tan\frac{\pi}{m}+s
$$
the sides of our $m-$polygon differ  from the regular $m-$polygon at most at $O(\sqrt{s}).$

Without loss of generality we may assume that
the largest and the smallest central angles, say $\angle U_{k-1}OU_{k}  = \frac{2\pi}{m} +2\phi $ and $  \angle U_{k}OU_{k+1}= \frac{2\pi}{m} +2\varphi , \ \varphi>\psi>0$  on Fig. 3  are adjacent.  Then the perimeter of the broken line $U_{k-1}V_{k-1} U_{k}V_{k}U_{k}$ is equal to $2\tan(\frac{\pi}{m} +2\varphi) +2 \tan(\frac{\pi}{m} -2\psi) .$ Making the symmetrization, we displace the point $U_k$ in the center of the arc $ U_{k-1}U_{k+1} $ and consider  the resulting new circumscribed $m-$polygon.

The perimeter of the  broken line $U_{k-1}V_{k-1} U_{k}V_{k}U_{k}$ is equal
 to $4\tan\left(\frac{\pi}{m} + \varphi-\psi \right), $ and the difference of two perimeters is $2\sin^2 (\varphi+\psi) /\cos(\frac{\pi}{m} +2\varphi)\cos(\frac{\pi}{m} -2\psi).$ It follows from this that $|\sin (\varphi+\psi) | < C_6 \sqrt{s},$ hence $\varphi+\psi < C_7\sqrt{s}.$
It follows that the inequalities  (\ref{differences}) and (\ref{alphas})  take place, possibly with other constants.

Now decompose the tangent of the sum of two angles
\begin{equation}
\label{summa}
\tan\left(\frac{\pi}{m}+\frac{\alpha_k-\alpha_{k-1}}{2}\right)=
\frac{\tan\frac{\pi}{m}+\tan\frac{\alpha_k-\alpha_{k-1}}{2}}{1-\tan\frac{\pi}{m}\tan\frac{\alpha_k-\alpha_{k-1}}{2}}.
\end{equation}
As the differences $|\alpha_k-\alpha_{k-1}|$ are small
 and as $\delta \to 0$
$$\tan \delta = \delta +\frac{1}{3}\delta^3 + O(\delta^4), \quad  \quad  \frac{1}{1-\delta}= 1+\delta+\delta^2+O(\delta^3),$$
we obtain, denoting temporarily $\delta:=\frac{\alpha_k-\alpha_{k-1}}{2},$ that
\begin{multline*}
\frac{\tan\frac{\pi}{m}+\tan\frac{\alpha_k-\alpha_{k-1}}{2}}{1-\tan\frac{\pi}{m}\tan\frac{\alpha_k-\alpha_{k-1}}{2}}
=\tan\frac{\pi}{m}+\delta\left(\tan^2\frac{\pi}{m}+1\right)+\delta^2\left(\tan^2\frac{\pi}{m}+1\right)\tan\frac{\pi}{m}+O(\delta^3). \\
\end{multline*}
By (\ref{a_k}) and (\ref{summa}), we get the following expression:
$$peri_m=2m\tan\frac{\pi}{m} + \frac{1}{2}(1+\tan^2\frac{\pi}{m})\tan\frac{\pi}{m}\sum\limits_{k=1}^{m-1}
(\alpha_k-\alpha_{k-1})^2+O\left(\sum\limits_{k=1}^{m-1}(\alpha_k-\alpha_{k-1})^3\right).$$

By means of  same reasoning as above  we
have
\begin{multline*}
P\left\{\sum\limits_{k=1}^{m-1}\alpha_k^2-\sum\limits_{k=1}^{m-1}\alpha_k\alpha_{k-1}\leq  \frac{s_5\cot\frac{\pi}{m}}{(1+\tan^2\frac{\pi}{m})}                                 \right\} \le\\
\le P\left\{peri_m \leq 2m \tan \frac{\pi}{m} +s,
\beta_1\leq\beta_2\leq\ldots\beta_{m-1}\right\} \le\\
\le P\left\{\sum\limits_{k=1}^{m-1}\alpha_k^2-\sum\limits_{k=1}^{m-1}\alpha_k\alpha_{k-1}\leq  \frac{s_6\cot\frac{\pi}{m}}{(1+\tan^2\frac{\pi}{m})}      \right\},
\end{multline*}
where $s_6>s_5$, and $s_5\to s, s_6 \to s$ as $s \to 0.$

Now denote for brevity
$$
w_k^2  = \frac{s_6\cot\frac{\pi}{m}} { \left(1+\tan^2 \frac{\pi}{m}\right) \left(1-\cos\left(\frac{\pi
k}{m}\right)\right)}, \ k=1,\dots, m-1.
$$

Using the orthogonal transform, we reduce the quadratic form under the sign of pro\-ba\-bility
to the form (\ref{orthogonal_transformation}), and obtain for small $s$:
\begin{multline*}
P\{peri_m \leq 2m\tan\frac{\pi}{m}+s\} \le \\
\le \cfrac{\Gamma(m)}{(2\pi)^{m-1}} \int_{\left\{  \sum\limits_{k=1}^{m-1}\left(1-\cos\frac{\pi
k}{m}\right)y_k^2   \le   \frac{s_{6} \cot \frac{\pi}{m} }{ \left(1+\tan^2 \frac{\pi}{m}\right) }    \right\}} \ dy_1\ldots
dy_{m-1} =\\
=\cfrac{\Gamma(m)}{(2\pi)^{m-1}} \int_{ \{\sum\limits_{k=1}^{m-1} \left(y_k/w_k \right)^2 \le1 \} }\ dy_1\ldots
dy_{m-1} = \frac{\Gamma(m)}{(2\pi)^{m-1}}b_{m-1}
\Pi^{''}_{m-1}(s_6),
\end{multline*}
where $b_{m-1}$ is from ~(\ref{volume_b_m}), and the
product of semi-axes of the new $(m-1)$-dimensional ellip\-so\-id  equals
$$ \Pi^{''}_{m-1}(s_6)=\prod\limits_{k=1}^{m-1}w_k =
s_6^{\frac{m-1}{2}}\left(\frac{\cot\frac{\pi}{m}}{2\left(1+\tan^2\frac{\pi}{m}\right)}\right)^{\frac{m-1}{2}}\frac{2^{m-1}}{\sqrt{m}}.$$

Thus, we obtain the following  estimate from above:
\begin{multline*}
s^{-\frac{m-1}{2}}P\{peri_m\geq 2m \tan\frac{\pi}{m}+s\} \le\\  \le
(s_6/s)^{\frac{m-1}{2}}\cfrac{\Gamma(m)}{(2\pi)^{m-1}}\left(\cfrac{\cot\frac{\pi}{m}}
{2\left(1+\tan^2\frac{\pi}{m}\right)}\right)^{\frac{m-1}{2}}\cfrac{2^{m-1}}
{\sqrt{m}}\cfrac{(\pi)^{\frac{m-1}{2}}}{\Gamma\left(\frac{m+1}{2}\right)}=\\
=(s_6/s)^{\frac{m-1}{2}}\cfrac{\Gamma(m)}{ \sqrt{m} \Gamma \left(\frac{m+1}{2}\right)}
\left(\frac{\cot\frac{\pi}{m}}{2\pi\left(1+\tan^2\frac{\pi}{m}\right)}\right)^{\frac{m-1}{2}}.
\end{multline*}
Taking the limit as ${s} \to 0,$ we obtain the upper estimate for the probability of interest. The lower estimate can be obtained analogously. The conclusion of Lemma 3.1 follows. \hfill $\square$

\textit{Proof of Theorem 3}. The proof is carried out similarly to the proof
of Theorem 1 and Theorem 2 by means of the transformation
$$z_n(t)=2m\tan\frac{\pi}{m}+tn^{-\frac{2m}{m-1}}.$$

In this case we verify the condition (\ref{corollary_silv_4}) for each $r\in \{1,\ldots,m-1\}$:
\begin{center}
$\lim\limits_{n\rightarrow\infty}n^{2m-r}P\{peri_n^{m;
1,\ldots,m}<z_n(t), peri_n^{m; m+1-r,\ldots,2m-r}<z_n(t)\}=0,$
\end{center}
which can be proved analogously to Lemma 2.2.

Basing on the conclusion of Silverman and Brown's Theorem (\ref{corollary_silv_3}) and considering the estimate (\ref{Lao_Mayer_Th})
 from Lao-Mayer Theorem, which can be reformulated in terms of $U$-min-statistic by replacement $h$ on $-h$,
i.e. $H_n=-\min(-h_n)$, we have:
$$\lim\limits_{n\rightarrow\infty}P\left\{peri_n^m>2m\tan\frac{\pi}{m} +tn^{-\frac{2m}{m-1}}\right\}=
\exp\left\{-\frac{t^{\frac{m-1}{2}}}{K_{3m}}\right\}.$$

Thus for each $t>0$
$$\lim\limits_{n\rightarrow\infty}P\left\{n^{\frac{2m}{m-1}}(peri_n^m
-2m\sin\frac{\pi}{m})\leq t\right\}=
1-\exp\left\{-\frac{t^{\frac{m-1}{2}}}{K_{3m}}\right\},$$ where
$$K_{3m}=m^{\frac{3}{2}}\left(2\pi\left(1+\tan^2\frac{\pi}{m}\right)\tan\frac{\pi}{m}\right)^{\frac{m-1}{2}}
\Gamma\left(\frac{m+1}{2}\right).$$
\hfill$\square$\medskip

The next result is a corollary of Theorem
\ref{theorem_3}.

\subsection{Area of circumscribed polygon}

{\bf Theorem 4}. {\it
Let $U_1, U_2, \ldots$ be independent and uniformly distributed points on the circum\-fe\-rence of
unit radius $\mathbb{S}$. For each set of the points $U_{i_1},\ldots,U_{i_m}$ on the circle
draw tangents at these points, denoting intersection of the tangents at the points $U_{i_k}$
and $U_{i_{k+1}}$  by $V_{i_k},  k=1,\ldots,m,$ correspondingly. Let
$area_n^m=\min\limits_{1\leq i_1<\ldots<i_m\leq n}
area_m(V_{i_1},\ldots,V_{i_m})$ be minimal of areas among all circumscribed $m$-polygons, formed
by points $V_{i_1},\ldots,V_{i_m}$. \ Then for each $t>0$
$$\lim\limits_{n\rightarrow\infty}P\left\{n^{\frac{2m}{m-1}}
(area_n^m-m\tan\frac{\pi}{m} )\leq
t\right\}=1-\exp\left\{-\frac{t^{\frac{m-1}{2}}}{K_{4m}}\right\},$$
where
$$K_{4m}=m^{\frac{3}{2}}\left(\pi\left(1+\tan^2\frac{\pi}{m}\right)\tan\frac{\pi}{m}\right)^{\frac{m-1}{2}}
\Gamma\left(\frac{m+1}{2}\right).$$ The rate of convergence is} $O\left(n^{-\frac{1}{ m-1 }}\right)$.

\textit{Proof of Theorem 4.} Note that  the minimal area and the minimal  perimeter of circumscribed $m$-polygon
in the case of unit circle  are related by the equality $area_n^{m}=\frac{1}{2}peri_n^{m}.$
Hence the following chain of equalities follows from Theorem \ref{theorem_3}:
\begin{multline*}
\lim\limits_{n\rightarrow\infty}P\left\{n^{\frac{2m}{m-1}}(peri_n^m
-2m\tan\frac{\pi}{m})\leq
t\right\}=\lim\limits_{n\rightarrow\infty}P\left\{n^{\frac{2m}{m-1}}(2area_n^m
-2m\tan\frac{\pi}{m})\leq t\right\}=\\
=\lim\limits_{n\rightarrow\infty}P\left\{n^{\frac{2m}{m-1}}(area_n^m
-m\tan\frac{\pi}{m})\leq t/2 \right\}.
\end{multline*}

On the other hand, by Theorem \ref{theorem_3} the right-hand side of the last equality is
$$1-\exp\left\{-\frac{t^{\frac{m-1}{2}}}{K_{3m}}\right\}=1-\exp\left\{-\frac{(t/2)^{\frac{m-1}{2}}}{2^{-\frac{m-1}{2}}K_{3m}}\right\}=
1-\exp\left\{-\frac{(t/2)^{\frac{m-1}{2}}}{K_{4m}}\right\},$$
where $K_{4m}=2^{-\frac{m-1}{2}}K_{3m},$ whence we get the conclusion of Theorem 4. \hfill$\square$

\section{Asymptotic behavior of constants}

In this section we investigate the asymptotics of the limit constants $K_{im}, i=\overline{1,4},$ as $m \to \infty,$
using the well-known Stirling's formula:  as $t \rightarrow +\infty$
$$\Gamma(t+1)\sim \sqrt{2\pi t}\left(\frac{t}{e}\right)^t, \,  t \rightarrow +\infty.$$

The constants from previous sections are as follows:
$$
\begin{array}{ll}
&K_{1m}=m^{\frac{3}{2}}\left(\pi\sin\frac{\pi}{m}\right)^{\frac{m-1}{2}}
\Gamma\left(\frac{m+1}{2}\right),\\
&K_{2m}=m^{\frac{3}{2}}\left(\pi\sin\frac{2\pi}{m}\right)^{\frac{m-1}{2}}
\Gamma\left(\frac{m+1}{2}\right)=\left(2\cos\frac{\pi}{m}\right)^{\frac{m-1}{2}}K_{1m},\\
&K_{3m}=m^{\frac{3}{2}}\left(2\pi\left(1+\tan^2\frac{\pi}{m}\right)\tan\frac{\pi}{m}\right)^{\frac{m-1}{2}}\Gamma\left(\frac{m+1}{2}\right)= 2^{\frac{m-1}{2}} K_{4m},\\
&K_{4m}=m^{\frac{3}{2}}\left(\pi\left(1+\tan^2\frac{\pi}{m}\right)\tan\frac{\pi}{m}\right)^{\frac{m-1}{2}}
\Gamma\left(\frac{m+1}{2}\right).
\end{array}
$$

Applying Stirling's formula, we get
$$
\Gamma\left(\frac{m+1}{2}\right)\sim\frac{\sqrt{2\pi\left(\frac{m-1}{2}\right)}
\left(\frac{m-1}{2}\right)^{\frac{m-1}{2}}}{e^{\frac{m-1}{2}}}=
\frac{\sqrt{\pi}(m-1)^\frac{m}{2}}{(2e)^{\frac{m-1}{2}}}\sim
\frac{\sqrt{\pi}m^{\frac{m}{2}}}{2^{\frac{m-1}{2}}e^{\frac{m}{2}}}, \, m \to \infty.
$$
The following relations clearly hold as $m\rightarrow\infty$:
$$\sin\frac{\pi}{m}\sim\frac{\pi}{m},
\quad\tan\frac{\pi}{m}\sim\frac{\pi}{m},\quad
1+\tan^2\frac{\pi}{m}\sim 1.$$
From here we easily deduce as $m\rightarrow\infty$:
$$
\begin{array}{ll}
&K_{1m}=m^{\frac{3}{2}}\left(\pi\sin\frac{\pi}{m}\right)^{\frac{m-1}{2}}
\Gamma\left(\frac{m+1}{2}\right)\sim
m^{\frac{3}{2}}\left(\pi\cdot\frac{\pi}{m}\right)^{\frac{m-1}{2}}
\frac{\sqrt{\pi}m^{\frac{m}{2}}}{2^{\frac{m-1}{2}}e^{\frac{m}{2}}}
=\frac{\pi^{m-\frac{1}{2}}m^2}{2^{\frac{m-1}{2}}e^{\frac{m}{2}}}:=\widetilde{K_{1m}},\\
&K_{2m}=\left(2\cos\frac{\pi}{m}\right)^{\frac{m-1}{2}}K_{1m}=m^{\frac{3}{2}}\left(\pi\sin\frac{2\pi}{m}\right)^{\frac{m-1}{2}}
\Gamma\left(\frac{m+1}{2}\right)\sim
\frac{\pi^{m-\frac{1}{2}}m^2}{e^{\frac{m}{2}}}:=\widetilde{K_{2m}}=2^{\frac{m-1}{2}}\widetilde{K_{1m}},\\
&K_{3m}=m^{\frac{3}{2}}\left(2\pi\left(1+\tan^2\frac{\pi}{m}\right)\tan\frac{\pi}{m}\right)^{\frac{m-1}{2}}
\Gamma\left(\frac{m+1}{2}\right)\sim
\frac{\pi^{m-\frac{1}{2}}m^2}{e^{\frac{m}{2}}}:=\widetilde{K_{3m}}=\widetilde{K_{2m}},\\
&K_{4m}=2^{-\frac{m-1}{2}}K_{3m}\sim \frac{\pi^{m-\frac{1}{2}}m^2}{2^{\frac{m-1}{2}}e^{\frac{m}{2}}}:=\widetilde{K_{4m}}=\widetilde{K_{1m}}.
\end{array}
$$
In total:

$$\widetilde{K_{2m}}=\widetilde{K_{3m}}=2^{\frac{m-1}{2}}\widetilde{K_{1m}}=2^{\frac{m-1}{2}}\widetilde{K_{4m}}=
\pi^{m-\frac{1}{2}} m^2 e^{\frac{m}{2}}.$$

Hence the asymptotic analysis shows the coincidence of the asymptotic constants
for the maximal area of inscribed random polygon and
for the minimal perimeter of circumscribed random polygon $\widetilde{K_{2m}}=\widetilde{K_{3m}},$
and also their coincidence for the maximal perimeter of inscribed random polygon and
for the minimal area of circumscribed random polygon. This seems to be the unexpected and curious observation. It would be interesting to understand if this follows from some general fact.

\section{Acknowledgements}

The research of the authors was supported by  the Program for Supporting Leading Scientific Schools (grant NSh-1216.2012.1). The authors are thankful to Dr. M. May\-er and  Dr.  Wei Lao  for the comments on some obscure points in their paper \cite{lao}.

\renewcommand{\refname}{References}

\center{Department of Mathematics and Mechanics, Saint-Petersburg State University, Russia}

\end{document}